\documentclass[leqno,12pt]{article}
\usepackage{latexsym}
\usepackage{graphicx}
\usepackage{amsmath}
\usepackage{amssymb}
\usepackage{mathtools}
\usepackage{cite}
\usepackage{color}
\usepackage{subfig}
\usepackage[margin=1.3in, top=1.3in, bottom=1.3in]{geometry}
\usepackage{algorithm}
\usepackage{algorithmic}
\usepackage{hyperref}

\newcommand{\nc}{\newcommand}
\nc{\nt}{\newtheorem}
\nt{thm}{Theorem}[section]
\nt{cor}[thm]{Corollary}
\nt{prop}[thm]{Proposition}
\nt{lem}[thm]{Lemma}
\nt{defn}[thm]{Definition}
\nt{rem}[thm]{Remark}
\nt{exa}[thm]{Example}
\nt{ass}[thm]{Assumption}
\nt{alg}[thm]{Algorithm}
\nt{conj}[thm]{Conjecture}
\nt{claim}[thm]{Claim}
\nt{oracle}[thm]{Oracle}
\nt{que}[thm]{Question}
\nc{\ip}[2]{\mbox{$\langle #1,#2 \rangle$}}
\nc{\pf}{\noindent{\bf Proof\ \ }}
\nc{\finpf}{\hfill{$\Box$}\linespace}
\nc{\linespace}{\vspace
{\baselineskip} \noindent}
\nc{\R}{{\mathbf R}}
\nc{\N}{{\mathbf N}}
\nc{\X}{{\mathbf X}}
\nc{\Y}{{\mathbf Y}}
\nc{\E}{{\mathbf E}}
\nc{\B}{{\mathbf B}}
\nc{\Sn}{{\mathbf S}}
\nc{\Hn}{{\mathbf H}}
\nc{\oR}{\overline{\R}}
\nc{\M}{\mathcal M}
\nc{\e}{\epsilon}

\nc{\Rn}{{\mathbf R}^n}
\nc{\inT}{\mbox{\rm int}\,}
\nc{\cl}{\mbox{\rm cl}\,}

\def\tto{\;{\lower 1pt \hbox{$\rightarrow$}}\kern -12pt
           \hbox{\raise 2.8pt \hbox{$\rightarrow$}}\;}
\newenvironment{myequation}{\setcounter{equation}{\value{thm}}
   \begin{equation}}{\addtocounter{thm}{1}\end{equation}}

\nc{\bmye}{\begin{myequation}}
\nc{\emye}{\end{myequation}}

\begin{document}
\title{
Strong growth and Goldstein subgradients in piecewise smooth optimization
}
\author{
Adrian S. Lewis
\thanks{ORIE, Cornell University, Ithaca, NY.
\texttt{people.orie.cornell.edu/aslewis} 
\hspace{2cm} \mbox{~}
Research supported in part by National Science Foundation Grant DMS-2405685.}
\and
Fahaar M. Pirani
\thanks{ORIE, Cornell University, Ithaca, NY.
\texttt{fp268@cornell.edu} }
}
\date{\today}
\maketitle

\begin{abstract}
We explore a new growth condition for Lipschitz functions around local minimizers, recently proposed in the context of Goldstein-subgradient-based algorithms and their behavior in practice.  On generic examples, these algorithms are often observed to converge approximately linearly.
We focus on objectives that are piecewise twice continuously differentiable.  In that case, the new growth condition holds when the objective is, in addition, strongly convex.  In the nonconvex case, we prove that quadratic growth in conjunction with a regularity property for the active gradients suffices.  Computational experiments illustrate the importance of the assumptions.
\end{abstract}
\medskip

\noindent{\bf Key words:} Lipschitz minimization, linear convergence, Goldstein \\
subgradients, quadratic growth, piecewise smooth 
\medskip

\noindent{\bf AMS Subject Classification:} 90C56, 49J52, 65Y20

\section{Introduction: growth and linear convergence}
A core concern for optimizers is the design of fast algorithms for minimizing  continuous objective functions $f$ over  Euclidean space $\X$.  We focus here on a locally Lipschitz function 
$f \colon \X \to \R$ with a strict local minimizer $\bar x \in \X$.  A classical goal is local {\em linear convergence}, meaning that the algorithm generates a sequence of iterates $x_k \to \bar x$ in $\X$, achieving an objective error $f(x_k) - f(\bar x)$ less than some given accuracy $\epsilon > 0$ within 
$k = \mbox{O}\big(\log(\frac{1}{\epsilon})\big)$ iterations.

Common in standard studies of linear convergence (surveyed in \cite{error_bound_d_lewis}) is  {\em quadratic growth}: 
\bmye \label{quadratic}
\mbox{
For some $\kappa > 0$, all $x \in \X$ near $\bar x$ satisfy 
$f(x)-f(\bar x) \ge \kappa |x-\bar x|^2$.
}
\emye
Closely related in classical scenarios is the following {\em linear subgradient growth} condition.  We denote the set of Clarke subgradients of $f$ at a point $x \in \X$ by $\partial f(x)$.
\bmye \label{growth}
\mbox{
For some $\alpha > 0$, all $x \in \X$ near $\bar x$ and $y \in \partial f(x)$ satisfy 
$|y| \ge \alpha |x-\bar x|$.
}
\emye
For a local minimizer $\bar x$ of a ${\mathcal C}^{(2)}$-smooth function $f$, for example, these two conditions are both equivalent to the Hessian $\nabla^2 f(\bar x)$ being positive definite.  

In the particular case of a convex function $f$ with a minimizer $\bar x \in \X$,  the two growth conditions (\ref{quadratic}) and (\ref{growth}) are also equivalent \cite[Theorem 3.5]{fran_sub}, and both hold in particular if $f$ is strongly convex.  In the case of the proximal point method, for example, \cite{luque} surveys convergence rates under various conditions, including linear subgradient growth (\ref{growth}).  

Rather than classical methods for convex or smooth optimization, our focus here is on recent algorithms for general Lipschitz optimization, originating with the seminal work \cite{mit} and followed up in \cite{ddmit,grimmer-goldstein,jordan-deterministic,kong-lewis,tian-so,somit,cutkovsky}.  These methods, rather than using proximal operations or subgradients at the current iterate $x_k$, rely instead on {\em Goldstein subgradients\/}: convex combinations of gradients (or Clarke subgradients) at points near $x_k$. The basic {\em Interpolated Normalized Gradient Descent (INGD)} algorithm of  \cite{mit} takes steps of some constant size $\delta>0$, terminating with a Goldstein subgradient of norm less than some tolerance $\epsilon > 0$.  Experiments on generic objectives suggest that {\em INGD with restarts}, using geometrically decaying step sizes $\delta = \epsilon$ (equal to $10^{-r}$ for $r=1,2,\ldots$ in our examples), usually results in linear convergence.  The following example is typical.

\begin{figure}
    \centering
    \includegraphics[width=0.7\textwidth]{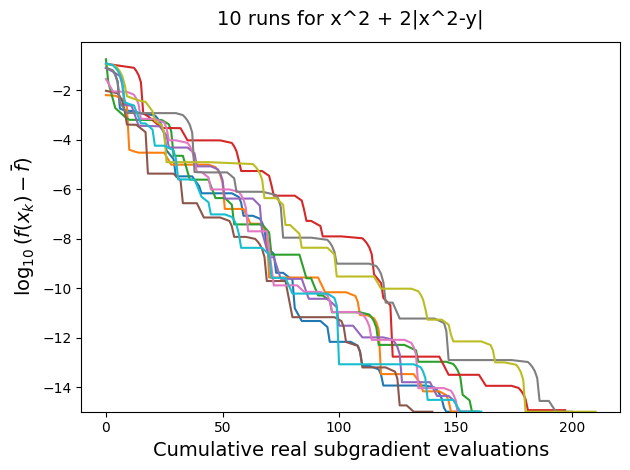}
    \caption{Ten randomly initialized runs of INGD with restarts on Example \ref{linear}, plotting objective error $\log_{10}(f - \min f)$ against subgradient calls.}
    \label{fig1}
\end{figure}

\begin{exa}[Linear convergence for a max function] \label{linear} \mbox{} \\
{\rm
The nonsmooth nonconvex function $f \colon \R^2 \to \R$ defined by
\[
f(x) ~=~ x_1^2 + 2|x_1^2 - x_2| \qquad (x \in \R^2)
\]
has unique global minimizer $\bar x = (0,0)$.  INGD with restarts \cite{mit} appears to converge approximately linearly relative to the number of subgradient evaluations:  see Figure~\ref{fig1}.
}
\end{exa}

\noindent
Our aim in this work is to explore under what conditions we might expect Goldstein subgradient-based algorithms like INGD with restarts to converge quickly, as in Example \ref{linear}. 

\subsection*{Contributions}
We study a novel enhancement to quadratic growth, a Goldstein-type growth condition introduced in \cite{goldstein-modulus}.  Motivated by several computational examples, we focus on piecewise twice continuously differentiable functions, and aim at sufficient conditions for the new growth condition.  We first prove, in Theorem \ref{strong}, that strong convexity suffices.  The nonconvex case, on the other hand, requires a more careful analysis.  In that case, in Theorem \ref{main}, we present a sufficient condition that combines the standard quadratic growth condition (\ref{quadratic}) with a simple generic regularity condition, easily verifiable in small instances.  Examples reveal the crucial nature of this regularity condition.  While the new growth condition has the potential to underpin fast convergence for nonsmooth nonconvex optimization \cite{goldstein-modulus}, whether it suffices specifically for the INGD method remains unknown.  Nonetheless, we hope that the theory and examples we present may prove illuminating. 

\subsection*{Organization}
In Section \ref{sec-flat}, we introduce the theme of our work, a reformulation of the growth condition of \cite{goldstein-modulus} in terms of ``Goldstein flatness''.  We present examples showing that strong convexity alone does not rule out such directions, and illustrating computationally that the INGD algorithm with restarts can converge slowly on strongly convex objectives.  Section \ref{sec-ruling} surveys previous results ruling out Goldstein flatness.  Section \ref{sec-piecewise} begins the discussion of piecewise smooth functions, proves that such functions, when strongly convex, cannot be Goldstein flat, and points out, via an example, the importance of strong convexity.  Section \ref{sec-linear} lays the linear-algebraic groundwork for a more careful study of piecewise smooth functions, leading into Section \ref{sec-orderings}, which studies the topology of the selection regions for piecewise smooth functions.  Finally, Section \ref{sec-gradient} presents the main result:  quadratic growth in conjunction with piecewise ${\mathcal C}^{(2)}$-smoothness and gradient regularity rules out Goldstein flatness.

\section{Goldstein flat directions} \label{sec-flat}
In experiments, Goldstein-subgradient-based methods like INGD with restarts very often seem to converge approximately linearly, just like in Example \ref{linear}.  However, a general rigorous analysis presents a challenge.  Remarkably, \cite{liwei} does develop a subtle version of the general approach that provably achieves {\em near linear} convergence, having complexity $\mbox{O}\big(\log^3(\frac{1}{\epsilon})\big)$, but that algorithm and its supporting conditions are quite involved.  Our aim here is a simple, broadly applicable growth condition.

To that end, we first note that the closely related conditions (\ref{quadratic}) and (\ref{growth}) ---
quadratic growth and linear subgradient growth --- are in general insufficient.  One underlying reason may be the surprising delicacy of linear subgradient growth.  Strongly convex functions satisfy both (\ref{quadratic}) and (\ref{growth}), and yet if we slightly modify condition~(\ref{growth}), replacing subgradients at the point $x$ by the expected gradient at a random point very near $x$, then this modified condition can fail even in the strongly convex case. The following example is an illustration.

\begin{figure}
    \centering
    \includegraphics[width=0.7\textwidth]{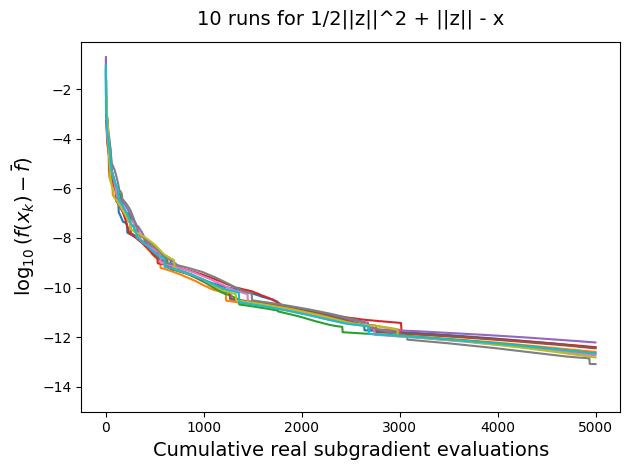}
    \caption{Ten randomly initialized runs of INGD with restarts on Example \ref{counterexample}, plotting objective error $\log_{10}(f - \min f)$ against subgradient calls.}
    \label{fig2}
\end{figure}

\begin{exa}[Strong convexity with small expected gradients] \label{counterexample} \mbox{} \\
{\rm
The strongly convex function $f \colon \R^2 \to \R$ defined by
\[
f(x) ~=~ \frac{1}{2}|x|^2 + |x| - x_1 \qquad (x \in \R^2)
\]
has a unique minimizer $\bar x = (0,0)$.  For each value $\tau \ge 0$, consider the point $x^\tau = (\tau,0)$.  With respect to the distance to the minimizer $|x^\tau - 0| = \tau$, the gradient norm  $|\nabla f(x^\tau)| = \tau$ grows linearly, as it must, $f$ being strongly convex.  However, consider a random variable $X^\tau$ uniformly distributed on the ball centered at $x^\tau$ with radius $(2\tau)^{3/2}$.  
Relative to its distance to the minimizer, this ball has small radius, but the expected gradient of $f$ at $X^\tau$ grows only quadratically:  ${\mathbb E}\big(\nabla f(X^\tau)\big) = \mbox{O}(\tau^2)$ as $\tau \downarrow 0$ (see Appendix \ref{A1}).  Interestingly, the convergence rate of INGD with restarts now appears much slower than linear:  see Figure \ref{fig2}.  
}
\end{exa}

Fast Goldstein-subgradient-based methods inspired by \cite{mit} were studied in \cite{goldstein-modulus}, under the assumption that the objective function $f$ has no ``Goldstein flat directions'' at the local minimizer $\bar x$, in the following sense. We denote the closed ball with radius $\e$ and center $x \in \X$ by $B_\e(x)$.

\begin{defn} \label{def-flat}
{\rm
A {\em Goldstein flat direction} for a locally Lipschitz function $f \colon \X \to \R$ at a point $\bar x \in \X$ is a unit vector $d \in \X$ for which there exist sequences of scalars $\tau_r > 0$ for $r=1,2,\ldots$ satisfying $\tau_r \downarrow 0$ as $r \to \infty$, and $\e_r = \mbox{o}(\tau_r)$, such that for each $r$, some convex combination $y_r$ of gradients of $f$ at points in the ball $B_{\e_r}(\bar x + \tau_r d)$ satisfies $|y_r| < \e_r$.  The function $f$ is {\em Goldstein flat} at $\bar x$ if it has a Goldstein flat direction.
}
\end{defn}

\noindent
In the terminology of \cite{goldstein-modulus}, the ``Goldstein modulus'' of $f$ grows linearly at 
$\bar x$ if and only if $f$ is not Goldstein flat at $\bar x$  (see Appendix \ref{goldstein}).  This novel growth condition is stronger than both quadratic growth and linear subgradient growth, as we observe next.

\begin{prop}
If a locally Lipschitz function is not Goldstein flat at a local minimizer, then both the quadratic growth condition (\ref{quadratic}) and the linear subgradient growth condition (\ref{growth}) hold.
\end{prop}

\pf
See Appendix \ref{goldstein}.
\finpf

\noindent
On the other hand, the converse of this result is false even for very well-structured functions, as a consequence of the following example.

\begin{exa}[Strong convexity with Goldstein flatness I] \label{strongly-convex3} \mbox{} \\
{\rm
The strongly convex function in Example \ref{counterexample} is semidefinite-representable \cite{mod_opt} and 
${\mathcal C}^{(2)}$-cone reducible \cite{Bon_Shap}, grows quadratically at its minimizer $\bar x = (0,0)$ and satisfies linear subgradient growth, but it has the Goldstein flat direction $d=(1,0)$.  Indeed, in the definition, we could choose $\tau_r = 1/r$, $\e_r = (2/r)^{3/2}$, and 
$y_r = {\mathbb E}\big(\nabla f(X^{\tau_r})\big)$.
}
\end{exa}

This example is atypical:  the minimizer $(0,0)$ for the function in Example \ref{counterexample}
is {\em degenerate}, in the following sense.  We denote the set of Clarke subgradients of a locally Lipschitz function $f \colon \X \to \R$ at a point $x \in X$ by $\partial f(x)$, and we say that $x$ is a {\em nondegenerate critical point} if $0 \in \mbox{ri}\big(\partial f(x)\big)$.  For a large class of objectives, critical points are generically nondegenerate \cite{gen_nondeg}: in particular, the sum of any convex or, more generally, lower-${\mathcal C}^{(2)}$ function with almost any linear perturbation has all critical points nondegenerate. 

\begin{figure}
    \centering
    \includegraphics[width=0.7\textwidth]{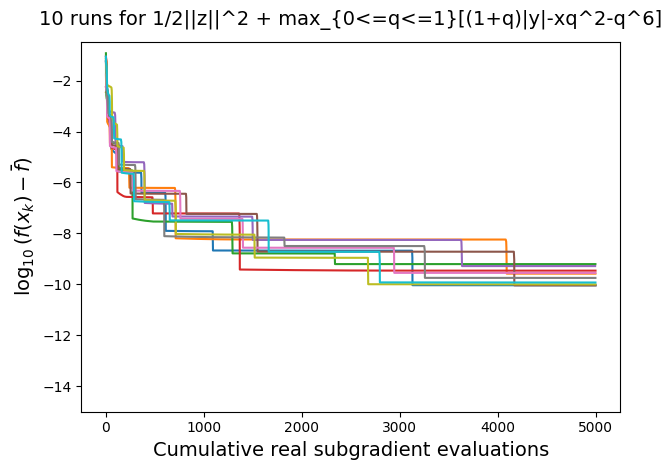}
    \caption{Ten randomly initialized runs of INGD with restarts on Example \ref{strongly-convex2}, plotting objective error $\log_{10}(f - \min f)$ against subgradient calls.}
    \label{figa}
\end{figure}

However, a more sophisticated example illustrates that the same behavior may occur even with nondegeneracy.
\begin{exa}[Strong convexity with Goldstein flatness II] \mbox{} \\
\label{strongly-convex2}
{\rm
The strongly convex function $f \colon \R^2 \to \R$ defined by
\[
f(x) ~=~ \frac{1}{2}|x|^2 ~+~ \max_{0 \le \theta \le 1} \{ (1+\theta)|x_2| - \theta^2 x_1 - \theta^6\} \qquad (x \in \R^2)
\]
is semi-algebraic, and is Goldstein flat at its minimizer, the nondegenerate critical point $\bar x = (0,0)$.  Indeed, a random variable $X^\tau$ uniformly distributed on the ball centered at the points 
$(\tau,0)$ with radius $(4\tau)^{3/2}$ satisfies ${\mathbb E}\big(\nabla f(X^\tau)\big) = \mbox{o}(\tau)$ as $\tau \downarrow 0$ (see Appendix \ref{E}).  Sublinear convergence of INGD with restarts is illustrated in Figure \ref{figa}.
}
\end{exa}

\section{Ruling out Goldstein flatness} \label{sec-ruling}
Ensuring that an objective is not Goldstein flat requires restrictive assumptions, as highlighted by Examples \ref{strongly-convex3} and  \ref{strongly-convex2}.  In this section we describe several sets of sufficient conditions. 

The first framework is a central focus of \cite{goldstein-modulus}.  We consider an objective $f$, not necessarily convex, that is a {\em ${\mathcal C}^{(2)}$ max function}, meaning a pointwise maximum of finitely-many ${\mathcal C}^{(2)}$-smooth functions.  Example \ref{linear} is an illustration.  

\begin{thm}[${\mathcal C}^{(2)}$ max functions {\rm \cite[Theorem 7.2]{goldstein-modulus}}] \mbox{} \\
Consider ${\mathcal C}^{(2)}$-smooth functions $f_i \colon \X \to \R$ for $i=1,2,\ldots,k$ and a point 
$\bar x$ where the values $f_i(\bar x)$ are all equal and the gradients $\nabla f_i(\bar x)$ are affinely independent.  Suppose that the function $f \colon \X \to \R$ is the pointwise maximum of the functions $f_i$, and that $\bar x$ is a nondegenerate critical point at which the quadratic growth condition~(\ref{quadratic}) holds.  Then $f$ is not Goldstein flat at $\bar x$.
\end{thm}

In passing, we touch on a second, very special framework, where the objective $f \colon \X \to \R$, again not necessarily convex, is {\em piecewise linear-quadratic}.  In other words, there are finitely-many polyhedra $X_i$ covering $\X$, and corresponding polynomials $q_i \colon \X \to \R$ of degree no larger than two, such that throughout each $X_i$ we have $f=q_i$.

\begin{prop} \label{plq}
A continuous, piecewise linear-quadratic function cannot be Goldstein flat at a strict local minimizer.
\end{prop}

\pf
See Appendix \ref{D}.
\finpf

A third framework known to rule out Goldstein flatness pertains to semi-algebraic or more generally stratifiable functions \cite{liwei}.  That approach is, however, most effective on objectives subject to generic linear perturbations \cite[Section 9]{goldstein-modulus}, whereas our interest is on particular instances.

\section{Piecewise smooth functions} \label{sec-piecewise}
To rule out Goldstein flatness, we must avoid apparently well-behaved objective functions like Examples \ref{strongly-convex3} and  \ref{strongly-convex2}.  Dispensing with traditional structural settings such as cone-representability or amenability, we focus instead on a simple and classical framework, highlighted in \cite{jongen-pallaschke}.  While this setting has some history from a variational analysis perspective, dating back to works such as \cite{kuntz-scholtes} and on through \cite{fac_pang}, its application in contemporary nonsmooth optimization specifically has been modest.  

\begin{defn} \label{pwdef}
{\rm
A function $f \colon \X \to \R$ is {\em piecewise ${\mathcal C}^{(2)}$} at a point $\bar x$ if there exist ${\mathcal C}^{(2)}$-smooth functions $f_i \colon \X \to \R$ for all $i$ in some finite index set $I$ such that for all points $x \in \X$ near $\bar x$, the value $f(x)$ lies in  the set $\{ f_i(x) : i \in I \}$.
}
\end{defn}

Examples of piecewise ${\mathcal C}^{(2)}$ functions include ${\mathcal C}^{(2)}$ max functions, continuous piecewise linear-quadratic functions, and, more generally, continuous {\em fully amenable} functions \cite{VA}.  This latter class consists of functions locally decomposable as $g \circ F$, for a convex piecewise linear-quadratic function $g \colon \R^n \to \R$ and a ${\mathcal C}^{(2)}$-smooth map $F \colon \X \to \R^n$.

\begin{prop}
If a continuous function $f \colon \X \to \R$ is piecewise ${\mathcal C}^{(2)}$ at a point $\bar x \in \X$, then it is locally Lipschitz around $\bar x$.
\end{prop}

\pf
See \cite{hager,kuntz-scholtes}.
\finpf

\begin{lem}
If a vector $a \in \X$ and a finite set $C \subset \X$ have the property that, for each point $x \in \X$, the inner product $\ip{a}{x}$ lies in the set $\{ \ip{c}{x} : c \in C \}$, then $a \in C$.
\end{lem}

\pf
By elementary linear algebra, there is a point $\bar x \in \X$ such that the inner products 
$\ip{c}{\bar x}$ (for $c \in C$) are all distinct.  Since the function $\ip{a}{\cdot}$ is continuous, there exists a vector $c \in C$ such that $\ip{a}{x} = \ip{c}{x}$ for all points $x$ near $\bar x$.  This implies $a = c$.
\finpf

\begin{prop}[Gradient selection] \label{gradient-selection}
For a function $f$ as in Definition~\ref{pwdef}, if $f$ is differentiable at a point $x \in \X$, then $\nabla f(x) = \nabla f_i(x)$ for some $i \in I$ satisfying $f(x) = f_i(x)$.
\end{prop}

\pf
For any vector $z \in \X$, and any sequence $0 < \tau_r \downarrow 0$, we have
\[
\ip{\nabla f(x)}{z} ~=~ \lim_{r \to \infty} \frac{1}{\tau_r}\big( f(x+\tau_r z) - f(x)\big).
\]
After taking a subsequence, we can suppose there exists $i \in I$ such that the right-hand side is
\[
\lim_{r \to \infty} \frac{1}{\tau_r}\big( f_i(x+\tau_r z) - f(x)\big).
\]
By continuity of $f_i$, we deduce $f(x) = f_i(x)$, and the limit above is then $\ip{\nabla f_i(x)}{z}$.  We have proved
\[
\ip{\nabla f(x)}{z} ~\in~ \{ \ip{\nabla f_i(x)}{z} : f_i(x) = f(x) \} \qquad \mbox{for all}~ z \in \X.
\]
The result now follows by the preceding lemma.
\finpf

\begin{prop}[Strongly convex piecewise smooth functions]  \label{strongpw}
Consider a $\mu$-strongly convex function $f \colon \X \to \R$ that is piecewise 
${\mathcal C}^{(2)}$ at a point $\bar x \in \X$.  Then any unit vector $u \in \X$, scalars $0 < \tau_r \downarrow 0$, and points $x_r = \bar x + \tau_r u + \mbox{o}(\tau_r)$ (for $r = 1,2,\ldots$) at which $f$ is differentiable satisfy 
\[
f'(\bar x;u) ~\le~  \ip{\nabla f(x_r)}{u} -\mu\tau_r + \mbox{o}(\tau_r) \qquad \mbox{as}~ r \to \infty.
\]
\end{prop}

\pf
Suppose that the function $f$ satisfies the conditions in Definition \ref{pwdef}.
For each $i \in I$, define an open set
\[
X_i ~=~ \{x \in \X : f(x) \ne f_j(x)~ \mbox{for all}~j \ne i \}.
\]
We can assume, without loss of generality, that for each $i \in I$ we have $\bar x \in \mbox{cl}\, X_i$.  If not, then we could replace the index set $I$ by $I \setminus \{i\}$ in Definition \ref{pwdef}, and repeating the process eventually results in a new set $I$ satisfying the desired property.

For each $i \in I$, there now exists a sequence $z_r \in X_i$, for $r=1,2,\ldots$, such that $z_r \to \bar x$.  We deduce $f_i(\bar x) = f(\bar x)$ for all $i \in I$.  Furthermore, since $f = f_i$ around $z_r$, and $f$ is $\mu$-strongly convex, we deduce 
$\nabla^2 f_i(z_r) \succeq \mu I$, and hence, by continuity, $\nabla^2 f_i(\bar x) \succeq \mu I$.

If the result fails, then, after taking a subsequence, we see that there exists a constant $\mu' < \mu$, a sequence $0 < \tau_r \downarrow 0$, and  a sequence of points $x_r = \bar x + \tau_r u + \mbox{o}(\tau_r)$, for $r = 1,2,\ldots$, such that 
\[
\ip{\nabla f(x_r)}{u} ~\le~ f'(\bar x;u) + \mu'\tau_r \qquad \mbox{for all}~ r.
\]
For each $r$, Proposition \ref{gradient-selection} implies the existence of an index $i \in I$ such that 
$\nabla f(x_r) = \nabla f_i(x_r)$ and $f(x_r) = f_i(x_r)$.  After taking a subsequence, we can assume that there exists $i \in I$ such that, for all $r=1,2,\ldots$, we have $f(x_r) = f_i(x_r)$ and 
\[
\ip{\nabla f_i(x_r)}{u} ~\le~ f'(\bar x;u) + \mu'\tau_r \qquad \mbox{for all}~ r.
\]
Since $f$ is convex, we have
\[
f'(\bar x;u) ~=~ \lim_{r \to \infty} \frac{1}{\tau_r} \big( f(x_r) - f(\bar x) \big)
~=~ \lim_{r \to \infty} \frac{1}{\tau_r} \big( f_i(x_r) - f_i(\bar x) \big).
\]
Since
\[
f_i(x_r) - f_i(\bar x) ~=~ \ip{\nabla f_i(\bar x)}{x_r - \bar x} + \mbox{o}(|x_r - \bar x|)
~=~ \tau_r\ip{\nabla f_i(\bar x)}{u} + \mbox{o}(\tau_r),
\]
we deduce $f'(\bar x;u) = \ip{\nabla f_i(\bar x)}{u}$.
We now have 
\begin{eqnarray*}
\ip{\nabla f_i(\bar x)}{u} +  \mu'\tau_r 
&\ge&
\ip{\nabla f_i(x_r)}{u} \\
&=&
\ip{\nabla f_i(\bar x) + \nabla^2 f_i(\bar x)(x_r - \bar x) + \mbox{o}(|x_r - \bar x|)}{u} \\
&=&
\ip{\nabla f_i(\bar x) + \tau_r \nabla^2 f_i(\bar x)u + \mbox{o}(\tau_r)}{u} \\
&\ge&
\ip{\nabla f_i(\bar x)}{u} + \mu\tau_r + \mbox{o}(\tau_r),
\end{eqnarray*}
which is our desired contradiction.
\finpf

The following new result, albeit restrictive in its assumptions, plays a central role motivating our development.

\begin{thm} \label{strong}
A strongly convex piecewise ${\mathcal C}^{(2)}$ function cannot be Goldstein flat at its minimizer.
\end{thm}

\pf
With the assumptions of Proposition \ref{strongpw}, suppose that the point $\bar x$ minimizes the function $f$ and some unit vector $u \in \X$ is a Goldstein flat direction.  In that case the directional derivative $f'(\bar x;u)$ must be nonnegative, and there exist scalars $0<\tau_r \downarrow 0$ and 
$\epsilon_r = \mbox{o}(\tau_r)$ such that some convex combination $y_r$ of gradients of $f$ at points in 
the ball $B_{\epsilon_r}(\bar x + \tau_r u)$ satisfies $|y_r| < \epsilon_r$ for all $r=1,2,\ldots$.  Each such gradient is $\nabla f(x_r)$ for some point $x_r \in B_{\epsilon_r}(\bar x + \tau_r u)$,
so as $r \to \infty$ we have, by Proposition \ref{strongpw},
\[
\ip{\nabla f(x_r)}{u}
~\ge~
f'(\bar x;u) + \mu\tau_r + \mbox{o}(\tau_r) ~\ge~  \mu\tau_r + \mbox{o}(\tau_r)
\]
and hence
\[
\epsilon_r ~>~ |y_r| ~\ge~ \ip{y_r}{u} ~\ge~ \mu\tau_r + \mbox{o}(\tau_r).
\]
This contradiction completes the proof.
\finpf

Two strong assumptions underlie the result of Theorem \ref{strong}, ruling out Goldstein flatness:  strong convexity and piecewise smoothness.  In general, even in the nondegenerate case, both assumptions are necessary.  On the one hand, Examples \ref{strongly-convex3} and \ref{strongly-convex2} show that we cannot replace piecewise smoothness by, say, semidefinite representability. On the other hand, the proof also makes crucial use of strong convexity.   

Instead of strong convexity, we might explore a weaker growth condition, common in contemporary variational analysis:  {\em tilt stability} \cite{tilt_orig}.  A local minimizer $\bar x$ is {\em tilt stable} for a locally Lipschitz function $f$ if, as we perturb the objective linearly to $f + \ip{y}{\cdot}$, the dependence of the unique corresponding local minimizer near $\bar x$ on the small vector $y$ should be Lipschitz.  For functions $f$ that are fully amenable, or more generally both subdifferentially continuous and prox-regular, tilt stability is equivalent to strong metric regularity of the subdifferential, and implies quadratic growth \cite[Theorem 3.3]{tilt}.   However, even convex fully amenable  functions can be Goldstein flat at a tilt-stable minimizer, as the following example shows.

\begin{figure}
    \centering
    \includegraphics[width=0.7\textwidth]{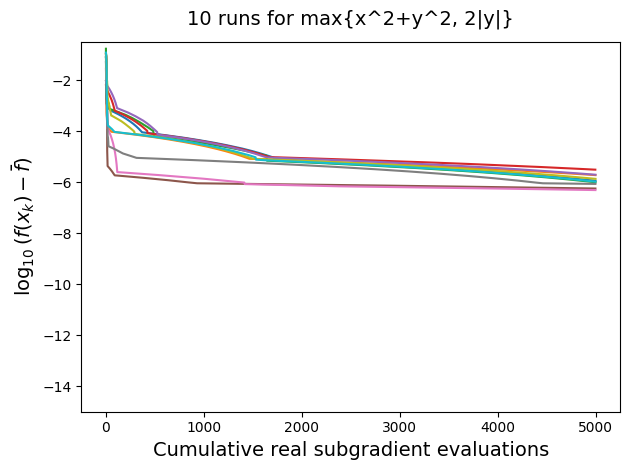}
    \caption{Ten randomly initialized runs of INGD with restarts on Example \ref{fully-amenable1}, plotting objective error $\log_{10}(f - \min f)$ against subgradient calls.}
    \label{figb}
\end{figure}

\begin{exa}[Tilt stability with Goldstein flatness] \label{fully-amenable1} \mbox{} \\
{\rm
Consider the convex function $f \colon \R^2 \to \R$ defined by
\[
f(x_1,x_2) ~=~ \max\{x_1^2+x_2^2,2|x_2|\Big\} \qquad \mbox{for}~ x \in \R^2.
\]
This function is fully amenable and semidefinite-representable, and its global minimizer $(0,0)$ is a nondegenerate critical point that is tilt stable.  However, $f$ is Goldstein flat at $(0,0)$.  Indeed, a random variable $X^\tau$ uniformly distributed on the ball centered at the points $(\tau,0)$ with radius $\tau^{3/2}$ satisfies ${\mathbb E}\big(\nabla f(X^\tau)\big) = \mbox{O}(\tau^{3/2})$ as $\tau \downarrow 0$ (see Appendix \ref{F}).  Slow convergence of INGD with restarts is illustrated in Figure \ref{figb}.
}
\end{exa}

\section{Regular lists of vectors} \label{sec-linear}
Our development concerns a continuous function $f \colon \X \to \R$ that is a {\em selection} of a set of smooth functions $f_i \colon \X \to \R$, indexed by $i$ in the set $[n] = \{1,2,\ldots,n\}$.  By this, we simply mean that, for each point $x \in \X$, there exists $i \in [n]$ such that $f(x) = f_i(x)$.  As $x$ varies, the selected index $j$ for which $f(x) = f_j(x)$ also varies.  However, if the values $f_i(x)$ for $i \in I$ are all distinct then, by continuity, $f = f_j$ around $x$.  The function $f$ is therefore well behaved over open subsets of $\X$ defined by strict orderings of functions values $f_i(x)$.  In the next section, we study those regions, via linearization.  For that purpose we need a regularity condition 
at the local minimizer $\bar x$ for $f$:  a linear-algebraic property of the list of gradients $\nabla f_i(\bar x)$.  We begin by developing that property. 

In the Euclidean space $\X$, we consider a list $G$ of $n$ vectors $g_i$, indexed by $i \in [n]$.  We view $G$ as an element of the space $\X^n$.  An {\em ordering} is just a permutation on $[n]$.  We call an ordering $\pi$ {\em strongly feasible} if there exists a vector $d \in \X$ satisfying
\[
\ip{g_{\pi(i+1)}}{d} ~<~ \ip{g_{\pi(i)}}{d} \qquad \mbox{for all}~ i \in [n-1],
\]
and {\em weakly feasible} if there exists a nonzero vector $d \in \X$ satisfying
\bmye \label{inequalities}
\ip{g_{\pi(i+1)}}{d} ~\le~ \ip{g_{\pi(i)}}{d} \qquad \mbox{for all}~ i \in [n-1].
\emye
An easy compactness argument shows the following result.

\begin{prop} \label{tool1}
For any ordering $\pi$, the set of lists $G \in \X^n$ for which $\pi$ is weakly feasible is closed.
\end{prop}

\pf
Consider a convergent sequence of lists $G^r \to G$ in $\X^n$.  For each $r=1,2,\ldots$, if the ordering $\pi$ is weakly feasible for $G^r$, then there exists a unit vector $d_r \in \X$ satisfying
\[
\ip{g^r_{\pi(i+1)}}{d_r} ~\le~ \ip{g^r_{\pi(i)}}{d_r} \qquad \mbox{for all}~ i \in [n-1].
\]
Any cluster point of the sequence $\{d_r\}$ then satisfies the inequalities (\ref{inequalities}), so $G$ is weakly feasible.
\finpf

If we define a corresponding list of vectors $G_\pi \in \X^{n-1}$ with elements
\[
(G_\pi)_i ~=~ g_{\pi(i+1)} - g_{\pi(i)} \qquad \mbox{for all}~ i \in [n-1], 
\]
then by separation we have the following easy result.  To simplify notation, we write $\mbox{\rm conv}\,G_\pi$ for the convex hull in $\X$ of the elements of $G_\pi$.

\begin{prop}[Feasibility via separation] \label{separation}
An ordering $\pi$ is strongly feasible for a list $G$ if and only if zero lies outside the polytope  
$\mbox{\rm conv}\,G_\pi$, and is weakly feasible if and only if zero lies outside the interior of 
$\mbox{\rm conv}\,G_\pi$.
\end{prop}

\pf
The result follows by a simple separating hyperplane argument.
\finpf

\begin{exa}[Strong and weak feasibility] \label{4-example}
{\rm
Consider the following list in $\R^2$:
\[
\left[ \begin{array}{cc} 2 \\ 1 \end{array} \right]
~,~
\left[ \begin{array}{cc} 0 \\ 1 \end{array} \right]
~,~
\left[ \begin{array}{cc} 1 \\ 0 \end{array} \right]
~,~
\left[ \begin{array}{cc} 1 \\ 3 \end{array} \right].
\]
There are 24 orderings.  Of these, exactly 12 are weakly feasible, and all of those are also strongly feasible.  Indeed, the 6 lines in $\R^2$ described by the equations $\ip{g_i-g_j}{x}=0$ for $1 \le i < j \le 4$ describe 12 cones in $\R^2$, each corresponding to a weakly feasible ordering.  In particular, the identity ordering is not weakly feasible.
}
\end{exa}

\begin{exa}[Weak versus strong feasibility] \label{3-example}
{\rm
Consider the following list in $\R^2$:
\[
\left[ \begin{array}{cc} 0 \\ 0 \end{array} \right]
~,~
\left[ \begin{array}{cc} 0 \\ 1 \end{array} \right]
~,~
\left[ \begin{array}{cc} 0 \\ -1 \end{array} \right].
\]
There are 6 orderings, all of which are weakly feasible.  However, the identity ordering is not strongly feasible.
}
\end{exa}

The following properties are easy to verify.

\begin{prop}[Affine independence]
If the vectors in a list are affinely independent, then all orderings are strongly feasible.
\end{prop}

\pf
If a list $G$ is affinely independent, then for any ordering $\pi$, the list $G_\pi$ is linearly independent, so the result follows quickly from the definition.
\finpf

\begin{prop}[Short lists]
For lists in $\X^n$ with $n \le \dim\X$, every ordering is weakly feasible.
\end{prop}

\pf
For any ordering $\pi$, the number of elements in the list $G_\pi$ is less than $\dim\X$, so their convex hull has empty interior.  Now we can apply Proposition \ref{separation}.
\finpf

The following idea is central to our development.

\begin{defn} \label{regular}
{\rm
A list $G \in \X^n$ is {\em regular} if, for every weakly feasible ordering $\pi$ on the set $[n]$, every subset of the list $G_\pi$ with size no more than $\dim\X$ is linearly independent.
}
\end{defn}
This idea is closely associated with the notion of ``full-spark'' matrices in linear algebra \cite{spark}.
Specifically, if we identify the space $\X$ with $\R^m$ and the list $G$ as an $m$-by-$n$ matrix, then $G$ is regular exactly when each of the $m$-by-$(n-1)$ matrices with columns given by the list $G_\pi$ for weakly feasible orderings $\pi$ is full-spark.  We remark that checking regularity is likely intractable for large instances (cf.\ \cite{tillmann-pfetsch}):  we are interested in the notion as a theoretical tool.

\begin{prop} \label{open}
The set of regular lists in $\X^n$ is open.
\end{prop}

\pf
Consider a sequence of irregular lists $G^r$ converging to a list $G$.  After taking a subsequence, we can suppose that there is an ordering $\pi$ that is weakly feasible for each $G^r$ and some fixed sublist from the list $G^r_\pi$ with size no more than $\dim\X$ that is linearly dependent.  In that case, $\pi$ is also weakly feasible for $G$, by Proposition~\ref{tool1}, and the corresponding sublist of the list $G_\pi$ is linearly dependent, so $G$ is also irregular.
\finpf

\begin{prop} \label{proper}
Consider any regular list $G \in \X^n$ and any weakly feasible ordering $\pi$.  Then no proper subspace of $\X$ contains a linearly dependent subset of the list $G_\pi$.
\end{prop}

\pf
By assumption, any linearly dependent subset $S$ of $G_\pi$ must have size exceeding $\dim\X$, and hence contains a subset of size $\dim\X$, which by regularity must be a basis for $\X$.  Consequently, $S$ is a spanning set, and hence is contained in no proper subspace.
\finpf

For the list in Example \ref{3-example}, the identity permutation is a weakly feasible ordering that is not strongly feasible.  However, this phenomenon is unusual, in the following sense.

\begin{thm}[Regularity and strong feasibility] \label{weak-strong}
\!\! For regular lists, all weakly feasible orderings are strongly feasible.
\end{thm}

\pf
Suppose that some ordering $\pi$ is weakly but not strongly feasible.  In that case zero lies on the boundary of the polytope $\mbox{conv}\,G_\pi$, and hence in one of its faces having dimension strictly less than $\dim\X$.  That face spans a proper subspace of $\X$, and contains a linearly dependent subset of $G_\pi$, contradicting Proposition \ref{proper}.
\finpf

The list in Example \ref{4-example} is regular.  On the other hand, for the list $G$ in Example~\ref{3-example}, the identity ordering $\iota$ is weakly feasible, and yet the matrix corresponding to $G_\iota$ is 
\[
\left[
\begin{array}{rr}
0 & 0 \\
1 & -2
\end{array}
\right],
\]
which is not full-spark, so $G$ is not regular.  However, regular lists are typical, in the following sense.

\begin{prop}[Regularity is generic]
The set of regular lists in $\X^n$ is open and dense.
\end{prop}

\pf
We have already observed, in Proposition \ref{open}, that the set of regular lists is open.  It remains to prove density.

Consider those lists of $(n-1)$ vectors in $\X$ for which every subset of size no larger than $\dim\X$ is linearly independent.  A standard argument concerning full-spark matrices \cite{blumensath} shows that such lists comprise a set $\Omega \subset \X^{n-1}$ that is open and dense.  For any ordering $\pi$, the linear map
from $\X^n$ to $\X^{n-1}$ defined by 
\[
(g_i)_{i=1}^n ~\mapsto~ (g_{\pi(i+1)} - g_{\pi(i)})_{i=1}^{n-1}
\]
is surjective and hence open.  Consequently, under this map, the inverse image of $\Omega$ is also open and dense.  The intersection of all such inverse images as $\pi$ ranges over the finitely-many orderings is therefore also open and dense.  The set of regular lists contains this set, so is also dense.
\finpf

\section{Orderings of smooth functions} \label{sec-orderings}
Equipped with our notion of regularity, we now consider the ingredients defining our piecewise smooth objective function.  In this section, we will only rely on ${\mathcal C}^{(1)}$-smoothness.  Specifically, we make the following assumption.

\begin{ass} \label{framework}
{\rm
For each $i \in [n]$, the function $f_i \colon \X \to \R$ is ${\mathcal C}^{(1)}$-smooth, and at the point $\bar x \in \X$, the values $f_i(\bar x)$ are all equal.   Furthermore, in $\X^n$, the list of gradients $\nabla f_i(\bar x)$ is regular.
}
\end{ass}

\begin{defn}
{\rm
For each ordering $\pi$, we define sets
\begin{eqnarray*}
X_\pi &=& \{ x \in \X : f_{\pi(i+1)}(x) \le f_{\pi(i)}(x) ~\mbox{for}~i \in [n-1]\} \\
\hat X_\pi &=& \{ x \in \X : f_{\pi(i+1)}(x) < f_{\pi(i)}(x) ~\mbox{for}~ i \in [n-1]\}.
\end{eqnarray*}
We call $\pi$ {\em feasible} if $\bar x \in \mbox{cl}\,\hat X_\pi$.
}
\end{defn}

\begin{prop} \label{equivalent}
If Assumption \ref{framework} holds, then for any ordering $\pi$, the following properties are equivalent.
\begin{enumerate}
\item[{\rm (i)}] $\pi$ is feasible.
\item[{\rm (ii)}] $\bar x \in \mbox{\rm cl}(X_\pi \setminus \{\bar x\})$.
\item[{\rm (iii)}] $\pi$ is weakly feasible for the list of gradients $\nabla f_i(\bar x)$.
\item[{\rm (iv)}] $\pi$ is strongly feasible for the list of gradients $\nabla f_i(\bar x)$.
\end{enumerate}
\end{prop}

\pf
Without loss of generality, suppose $\bar x = 0$.
Clearly property (i) implies property (ii) because $\hat X_\pi \subset X_\pi \setminus \{0\}$.  If (ii) holds, then there exists a sequence of nonzero vectors $x_r \in X_\pi$ converging to $0$.  After taking a subsequence, we can suppose that the normalized vectors $\frac{1}{|x_r|} x_r$ converge to some unit vector $d \in \X$.  We deduce $\ip{\nabla f_{\pi(i+1)}(0)}{d} \le \ip{\nabla f_{\pi(i)}(0)}{d}$ for all $i \in [n-1]$, so property (iii) holds.  Theorem \ref{weak-strong} ensures that (iii) implies property (iv).  Finally, if (iv) holds, then some vector $d \in \X$ satisfies $\ip{\nabla f_{\pi(i+1)}(0)}{d} < \ip{\nabla f_{\pi(i)}(0)}{d}$ for all $i \in [n-1]$, so $\tau d \in \hat X_\pi$ for all small $\tau > 0$, and hence (i) holds.
\finpf

\begin{prop}
If Assumption \ref{framework} holds, then 
\[
\bar x ~\in~ \mbox{\rm int}\Big(  \bigcup_{\mbox{\scriptsize\rm feasible}~\pi} X_\pi \Big).
\]
\end{prop}

\pf
If the result fails, then there exists a sequence $(x_r)$ disjoint from $X_\pi$ for each feasible ordering $\pi$ and converging to the point $\bar x$.  Clearly, each $x_r$ is nonzero, and must lie in a set $X_{\hat\pi}$ for some ordering $\hat\pi$.  After taking a subsequence, we can suppose that the entire sequence is contained in one set $X_{\hat\pi}$.  However, in that case $\hat\pi$ is feasible, by Proposition \ref{equivalent}.
\finpf

Given a Euclidean space $\Y$, a function $h \colon \Y \to \R$ is {\em approximately convex} at a point $v \in \Y$ if, for all $\epsilon > 0$, 
\[
h(\lambda y + (1-\lambda)z) ~\le~ \lambda h(y) + (1-\lambda)h(z) + \epsilon \lambda(1-\lambda)|y-z|
\]
for all points $y,z \in \Y$ near $v$ and all $\lambda \in [0,1]$.  Any ${\mathcal C}^{(1)}$-smooth function is approximately convex, as is any pointwise maximum of finitely-many ${\mathcal C}^{(1)}$-smooth functions, and, more generally, any function that is lower-${\mathcal C}^{(1)}$ at $v$:  see \cite{daniilidis-georgiev}.

We need the following simple tool.

\begin{lem} \label{lemma-strict}
Consider a Euclidean space $\Y$ and a function $h \colon \Y \to \R$ that is approximately convex at $0$ and satisfies $h(0) = 0$.  All nonzero vectors $y \in \Y$ near $0$, and scalars $\tau \ge h(y)$ and $\lambda \in [0,1)$ satisfy
\[
h(\lambda y) ~<~ \sqrt{\tau^2 + |y|^2(1-\lambda^2)}.
\]
\end{lem}

\pf
Approximate convexity implies
\[
h(\lambda y) ~=~ h(\lambda y + (1-\lambda)0) 
~\le~ \lambda h(y) + (1-\lambda)h(0) + \lambda(1-\lambda)|y-0|
~=~ \lambda h(y) + \lambda(1-\lambda)|y|.
\]
We can assume $\lambda > 0$, since otherwise the result is immediate.
Strict concavity then implies
\[
\sqrt{\tau^2 + |y|^2(1-\lambda^2)} ~>~ \lambda|\tau| + (1-\lambda)\sqrt{\tau^2 + |y|^2}
~\ge~ \lambda h(y) + (1-\lambda)|y|,
\]
and the result follows.
\finpf

\begin{prop}
\mbox{}~~Given a Euclidean space $\Y$, consider a continuous function \mbox{$h \colon \Y \to \R$} that is approximately convex at $0$ and satisfies $h(0) = 0$.  For any small radius $\rho > 0$, the sets
\begin{eqnarray*}
\hat X &=& \{ (y,\tau) \in \Y \times \R : h(y) < \tau,~ |y|^2 + \tau^2 < \rho^2 \}  \\
\overline X &=& \{ (y,\tau) \in \Y \times \R : h(y) \le \tau,~ |y|^2 + \tau^2 \le \rho^2 \}
\end{eqnarray*}
are nonempty and connected, and satisfy $\overline X = \mbox{\rm cl}\,\hat X$ and $\hat X = \mbox{\rm int}\,\overline X$.
\end{prop}

\pf
Clearly the set $\overline X$ is closed and $\hat X \subset \overline X$, so $\mbox{cl}\,\hat X \subset \overline X$.  On the other hand, the set $\hat X$ is open and $\hat X \subset \overline X$, so $\hat X \subset \mbox{int}\, \overline X$.

Consider any point $(y,\tau) \in \hat X$.  The path
\[
 \lambda ~\mapsto~ (\lambda y, \sqrt{\tau^2 + |y|^2(1-\lambda^2)})
 \]
 as $\lambda$ increases from $0$ to $1$ connects a point on the open line segment 
 $\{0\} \times (0,\rho) \subset \hat X$ to the point $(y,|\tau|)$, and lies in $\hat X$, by the preceding lemma.  The line segment (possibly trivial) between the points $(y,\tau)$ and $(y,|\tau|)$ also lies in 
 $\hat X$.  This proves that $\hat X$ is connected.
 
 We next show that any point $(y,\tau) \in \overline X$ is a limit of points in $\hat X$.  If $\tau < 0$, then we note $(y,\tau') \in \hat X$ for all $\tau' \in (\tau, 0)$.  If $y=0$, then either $\tau = 0$ in which case we note $(0,\tau') \in \hat X$ for all small $\tau' > 0$, or $\tau > 0$, in which case we note 
 $(0,\tau') \in \hat X$ for all $\tau' \in (0,\tau)$.  Finally, if $y \ne 0$ and $\tau > 0$, then, using Lemma \ref{lemma-strict}, we can choose scalars $0 < \theta_r \to \tau$ as $r \to \infty$, satisfying
\[
h( \big(\cos(1/r)\big))y) ~<~ \theta_r ~<~ \sqrt{\tau^2 + |y|^2\sin^2(1/r)} 
\]
and then  we note that the points $(\big(\cos(1/r)\big)y, \theta_r) \in \hat X$ converge to the point $(y,\tau)$.  We deduce that the set $\overline X$ is the closure of $\hat X$, and hence is also connected.
 
 It remains to prove $\mbox{int}\, \overline X \subset \hat X$, so consider any point 
 $(y,\tau) \in \mbox{int}\, \overline X$.  For all $\tau'$ near $\tau$ we have  
 $(y,\tau') \in \overline X$ and hence $h(y) \le \tau'$, so $h(y) < \tau$.  Furthermore, for all $\lambda$ near $1$ we have $\lambda(y,\tau) \in \overline X$ and hence $\lambda^2(|y|^2 + \tau^2) \le \rho^2$, so
 $|y|^2 + \tau^2 < \rho^2$, completing the proof.
 \finpf

\begin{prop}[Local topology of feasible regions] \label{tool2}
Consider a  ${\mathcal C}^{(1)}$-smooth function $F \colon \X \to \R^m$ satisfying $F(0) = 0$.  Suppose that there exists a  vector $d \in \X$ satisfying $DF(0)d < 0$.  Then for any small radius $\rho > 0$, the sets
\begin{eqnarray*}
\bar X &=& \{ x \in \X : F(x) \le 0,~ |x| \le \rho\} \\
\hat X &=& \{ x \in \X : F(x) < 0,~ |x| < \rho\}
\end{eqnarray*}
are nonempty and connected, and satisfy $\bar X = \mbox{\rm cl}\,\hat X$ and $\hat X = \mbox{\rm int}\,\bar X$.
\end{prop}

\pf
Without loss of generality, suppose $|d|=1$.  Denote the orthogonal complement of the vector $d$ by $\Y$.  Then we can write $\X = \Y \times \R$, identifying any pair $(y,\tau) \in \Y \times \R$ with the point 
$y+\tau d$.  By the implicit function theorem, there exist ${\mathcal C}^{(1)}$-smooth functions 
$h_i \colon \Y \to \R$ for all $i \in [m]$ satisfying
\begin{eqnarray*}
F_i(y,\tau) \le 0 & \Leftrightarrow & h_i(y) \le \tau \\
F_i(y,\tau) < 0 & \Leftrightarrow & h_i(y) < \tau
\end{eqnarray*} 
for all points $y \in \Y$ near $0$ and all small $\tau \in \R$.  The continuous function $h \colon \Y \to \R$ defined by $h(y) = \max_i h_i(y)$ is approximately convex at $0$, and for any small radius $\rho >0$ we have
\begin{eqnarray*}
\bar X &=& \{ (y,\tau) \in \Y \times \R : h(y) \le \tau,~ |y|^2 + \tau^2 \le \rho^2 \} \\
\hat X &=& \{ (y,\tau) \in \Y \times \R : h(y) < \tau,~ |y|^2 + \tau^2 < \rho^2 \}.
\end{eqnarray*}
The result now follows from the preceding Proposition.
\finpf

Using this tool, we deduce the following property.

\begin{prop} \label{intcl}
Suppose that Assumption \ref{framework} holds.
Fix any small radius $\rho > 0$.  For any feasible ordering $\pi$, consider the sets
\[
X_\pi^\rho ~=~ X_\pi \cap B_\rho(\bar x) \qquad \mbox{and} \qquad 
\hat X_\pi^\rho ~=~ \hat X_\pi \cap \mbox{\rm int}\, B_\rho(\bar x).
\]
The two sets have equal measure, satisfy $\hat X_\pi^\rho = \mbox{\rm int}\,X_\pi^\rho$ and $X_\pi^\rho = \mbox{\rm cl}\,\hat X_\pi^\rho$, and both are connected.
\end{prop}

\pf
Without loss of generality we can suppose $\bar x = 0$ and $f_i(\bar x) = 0$ for all $i \in [n]$.
With the exception of the equality of measure, the proof follows from Proposition~\ref{tool2}.  Each boundary point lies either on the sphere of radius $\rho$ or in one of the sets
\[
\{x \in \X : f_{\pi(i+1)}(x) = f_{\pi(i)}(x) \} \qquad \mbox{for}~ i \in [n-1].
\]
By regularity, $\nabla f_{\pi(i+1)}(0) \ne \nabla f_{\pi(i)}(0)$, so each such set is a proper submanifold of $\X$ around $0$, and hence has measure zero near $0$.
\finpf

We now consider a continuous selection of the functions $f_i$. 

\begin{prop}[Constant selection] \label{constant}
If Assumption \ref{framework} holds, then
for each feasible ordering $\pi$ there exists a unique index \mbox{$\phi(\pi) \in [n]$} such that \mbox{$f(x)=f_{\phi(\pi)}(x)$} for all points $x \in X_\pi$ near the point $\bar x$.  Furthermore, if $f$ is differentiable at such a point $x$, then $\nabla f(x) = \nabla f_{\phi(\pi)}(x)$.
\end{prop}

\pf
Choose a radius $\rho > 0$ sufficiently small to apply Proposition \ref{intcl}, and consider any feasible ordering $\pi$.  For each point $x \in \hat X^\rho_\pi$, there exists a unique index $i \in [n]$ satisfying $f(x) = f_i(x)$:  we denote this $i$ by $\psi(x)$.  By continuity, this selection function $\psi \colon \hat X_\pi \to [n]$ is locally constant throughout the connected open set $\hat X_\pi^\rho$, and hence is constant.  We denote the value by $\phi(\pi)$.  Applying Proposition~\ref{intcl} and continuity again shows $f = f_{\phi(\pi)}$ throughout the set $X_\pi^\rho$, as required.  For the final claim, since the ordering $\pi$ is strongly feasible for the list of gradients $\nabla f_i(\bar x)$, that remains the case for $\nabla f_i(x)$, so the set $X_\pi$ has a full-dimensional tangent cone at $x$.  Since $f$ is differentiable at $x$, and $f = f_{\phi(\pi)}$ throughout the set $X_\pi^\rho$, the result follows.
\finpf

\noindent
As a consequence of this result, we see that, around the point $\bar x$, the function $f$ is \mbox{${\mathcal C}^{(1)}$-smooth} throughout the dense, open, full-measure set 
\[
\Omega ~=~ \bigcup_{\mbox{\scriptsize feasible}~\pi} \hat X_\pi .
\]

\section{Gradient regularity and Goldstein flatness} \label{sec-gradient}
The function $f \colon \R^2 \to \R$ in Example \ref{fully-amenable1} is a continuous selection of the 
${\mathcal C}^{(2)}$-smooth functions
\[
f_1(x_1,x_2) ~=~ |x|^2, \qquad f_2(x_1,x_2) ~=~ 2x_2, \qquad f_3(x_1,x_2) ~=~ -2x_2.
\]
It grows quadratically at $0$, and yet is Goldstein flat there.  As we shall see, the missing ingredient for ruling out Goldstein flatness is regularity of the gradients $\nabla f_i(0)$, which fails in this example.  We first prove a simple tool.

\begin{lem} \label{tool}
If a ${\mathcal C}^{(2)}$-smooth function $h \colon \R \to \R$ satisfies $h(0)=0$ and $h(t) \ge \frac{\kappa}{2} t^2$ for all small $t \ge 0$, then $h'(t) \ge \frac{\kappa}{2} t$ for all small $t \ge 0$.
\end{lem}

\pf
Clearly $h'(0) \ge 0$.  If $h'(0) > 0$, the result follows by continuity of $h'$.  On the other hand, if $h'(0) = 0$, then the Taylor expansion shows $h''(0) \ge \kappa$, and the result then follows.
\finpf

We can now prove our main result.

\begin{thm}[Piecewise smoothness and Goldstein flatness]~~ \label{main}
Consider a continuous piecewise ${\mathcal C}^{(2)}$ function $f \colon \X \to \R$ that satisfies the quadratic growth property (\ref{quadratic}) at the point $\bar x \in \X$.  Suppose furthermore that the following regularity property holds:  $f$ is a continuous selection of ${\mathcal C}^{(2)}$-smooth functions $f_i \colon \X \to \R$, for $i \in [n]$, and the list of gradients $\nabla f_i(\bar x)$ for those $i \in [n]$ satisfying $f_i(\bar x) = f(\bar x)$ is regular.  In that case, $f$ is not Goldstein flat at $\bar x$.
\end{thm}

\pf
We lose no generality in supposing that $\bar x = 0$.  After discarding those functions $f_i$ not satisfying $f_i(0) = f(0)$, we can also suppose that $f_i(0) = 0$ for all $i \in [n]$.
Suppose that the property described in Definition \ref{def-flat} of Goldstein flatness holds.  

Consider the equivalence relation 
$\sim$ on the set $[n]$ defined, via the Goldstein flat direction $d \in \X$, by $i \sim j$ if and only if 
$\ip{\nabla f_i(0)}{d} = \ip{\nabla f_j(0)}{d}$.  We claim that the set
\[
\M ~=~ \{x \in \X : f_i(x)=f_j(x) ~\mbox{whenever}~ i \sim j\}
\]
is a ${\mathcal C}^{(2)}$-manifold around $0$.  To see this, choose any ordering $\pi$ such that 
\[
\ip{\nabla f_{\pi(i+1)}(0)}{d} ~\le~ \ip{\nabla f_{\pi(i)}(0)}{d} \qquad \mbox{for all}~ i \in [n-1].
\]
By definition, $\pi$ is weakly feasible.
Consider the set
\[
\M' ~=~ \{x \in \X : f_{\pi(i+1)}(x)=f_{\pi(i)}(x) ~\mbox{whenever}~ \pi(i+1) \sim \pi(i)\}
\]
Clearly $\M \subset \M'$.  On the other hand, consider any point $x \in \M'$, and suppose $i \sim j$. Define $i' = \pi^{-1}(i)$ and $j' = \pi^{-1}(j)$, and suppose, without loss of generality, $i' < j'$.
Since $\ip{\nabla f_{\pi(i')}(0)}{d} = \ip{\nabla f_{\pi(j')}(0)}{d}$, in fact we have
\[
\ip{\nabla f_{\pi(k+1)}(0)}{d} ~=~ \ip{\nabla f_{\pi(k)}(0)}{d},
\]
or in other words $\pi(k+1) \sim \pi(k)$, for $i' \le k < j'$.  Since $x \in \M'$ we deduce 
$f_{\pi(k+1)}(x) = f_{\pi(k)}(x)$ for $i' \le k < j'$, so $f_i(x) = f_j(x)$.  We have therefore proved $\M = \M'$.  The orthogonal complement of $d$ (a proper subspace of $\X$), contains the set of vectors 
\[
\{ \nabla f_{\pi(i+1)}(0) - \nabla f_{\pi(i)}(0) : \pi(i+1) \sim \pi(i) \}
\]
which Proposition \ref{proper} implies is linearly independent, by regularity.  Thus $\M$ is a manifold around~$0$.

Notice that the Goldstein flat direction $d$ lies in the tangent space to the manifold $\M$ at $0$.  Hence there exists a ${\mathcal C}^{(2)}$-smooth path $p \colon (-\delta,\delta) \to \M$, for some $\delta > 0$, such that $p(0) = 0$ and $p'(0) = d$.

Now consider any sequence of points $x_r \in \X$ at which $f$ is differentiable and satisfying $|x_r - \tau_r d| = \mbox{o}(\tau_r)$.  We claim
\bmye \label{key-claim}
\ip{\nabla f(x_r)}{p'(\tau_r)} ~\ge~ \frac{\kappa}{3} \tau_r \qquad \mbox{for all large $r$.}
\emye

By way of contradiction, suppose that property (\ref{key-claim}) fails.  After taking a subsequence, we can suppose there exists an ordering $\pi$, and a sequence of points $x_r \in X_\pi$  such that 
\[
\ip{\nabla f(x_r)}{p'(\tau_r)} ~<~ \frac{\kappa}{3} \tau_r \qquad \mbox{for all}~r.
\]
By Proposition \ref{constant}, we have $\nabla f(x_r) = \nabla f_{\phi(\pi)}(x_r)$ for all large $r$.
Using Proposition~\ref{intcl}, we can perturb each point $x_r$ slightly, and thereby suppose in addition that $x_r \in \hat X_\pi$ for all $r$.

Continuing our argument, we now claim
\bmye \label{claima}
p(\tau) \in X_\pi \qquad \mbox{for all small $\tau \ge 0$}.  
\emye
If not, then there exists a sequence $0 < \tau'_r \downarrow 0$ such that, for all $r=1,2,\ldots$ there exists an index $i \in [n-1]$ with 
\[
f_{\pi(i+1)}\big(p(\tau'_r)\big) > f_{\pi(i)}\big(p(\tau'_r)\big).
\]
After taking a subsequence, we can suppose that there is a fixed index $i$ for which this inequality holds for all $r$.  If $\pi(i+1) \sim \pi(i)$, then 
\[
f_{\pi(i+1)}\big(p(\tau)\big) = f_{\pi(i)}\big(p(\tau)\big) \qquad \mbox{for all small}~ \tau \ge 0,
\]
which is a contradiction, so the quantity
\[
\beta ~=~ \ip{\nabla f_{\pi(i+1)}(0) - \nabla f_{\pi(i)}(0)}{d}
\]
is nonzero.  If $\beta < 0$, then
\[
f_{\pi(i+1)}\big(p(\tau)\big) - f_{\pi(i)}\big(p(\tau)\big) ~<~ 0 \qquad \mbox{for all small}~ \tau \ge 0,
\]
which is again a contradiction.  We deduce $\beta > 0$.  We next observe
\[
f_{\pi(i+1)}\big(p(\tau)\big) - f_{\pi(i)}\big(p(\tau)\big) ~>~ \frac{\beta}{2} \tau \qquad \mbox{for all small}~ \tau > 0.
\]
Furthermore, we have $|x_r - p(\tau_r)| = \mbox{o}(\tau_r)$, giving the contradiction
\[
\frac{\beta}{2} \tau_r 
~<~ f_{\pi(i+1)}\big(p(\tau_r)\big) - f_{\pi(i)}\big(p(\tau_r)\big)
~=~ f_{\pi(i+1)}(x_r) - f_{\pi(i)}(x_r) + \mbox{o}(\tau_r)
~\le~ \mbox{o}(\tau_r)
\]
as $r \to \infty$.  We have thus proved property (\ref{claima}).

The ordering $\pi$ is feasible, so $f=f_{\phi(\pi)}$ throughout the set $X_\pi$ near $0$.  Since 
$|p(\tau)| = \tau + \mbox{o}(\tau)$ for all small $\tau \ge 0$, quadratic growth ensures
\[
f_{\phi(\pi)}\big(p(\tau)\big) ~=~ f\big(p(\tau)\big) ~\ge~ \kappa|p(\tau)|^2 ~\ge~ \frac{\kappa}{2} \tau^2
\]
for all small $\tau \ge 0$.

Applying Lemma \ref{tool} to the composite function $f_{\phi(\pi)} \circ p$ shows 
\[
\ip{\nabla f_{\phi(\pi)} \big(p(\tau)\big)}{p'(\tau)} ~\ge~ \frac{\kappa}{2} \tau
\]
for all small $\tau \ge 0$.  Since $|x_r - p(\tau_r)| = \mbox{o}(\tau_r)$, we deduce
\[
\ip{\nabla f(x_r)}{p'(\tau_r)} ~=~ 
\ip{\nabla f_{\phi(\pi)}(x_r)}{p'(\tau_r)} ~\ge~ \frac{\kappa}{3} \tau_r
\]
for all large $r$, contradicting our assumption that inequality (\ref{key-claim}) fails.

We have thus proved inequality (\ref{key-claim}).
Any convex combination $g_r$ of gradients of $f$ at points in balls of the form $B_{\epsilon_r}(\tau_r d)$ with radius $\epsilon_r = \mbox{o}(\tau_r)$  therefore satisfies
$\ip{g_r}{p'(\tau_r)} \ge \frac{\kappa}{3} \tau_r$, for large $r$.  Since the vectors $p'(\tau_r)$ converge to the unit vector $d$, we deduce \mbox{$|g_r| \ge \frac{\kappa}{4} \tau_r$,} for large $r$, contradicting Goldstein flatness.
\finpf


{\small
\parsep 0pt
\def\cprime{$'$} \def\cprime{$'$}

}

\appendix
\section{Appendix:  Proof of Example \ref{counterexample}} \label{A1}
\begin{lem}
If a random vector $X$ is uniformly distributed over the disk in $\R^2$ of radius $\rho \ge 0$ centered at the first unit vector $e_1$, then
\[
{\mathbb E}\frac{X_1}{|X|} ~=~ 1 - \frac{\rho^2}{8} + \mbox{\rm O}(\rho^4) \qquad \mbox{as}~ \rho \downarrow 0.
\]
\end{lem}

\pf
For radius $\rho \in [0,1)$ we have 
\[
\beta ~=~\int_{|x-e_1| \le \rho} \frac{x_1}{|x|}\, dx ~=~ \int\!\!\!\int_{u^2+v^2 \le \rho^2} f(u,v)\,du\,dv
\]
where the function $f$ defined, for points $(u,v) \in B_{\rho}(0,0)$, by
\[
f(u,v) ~=~ \frac{1+u}{\sqrt{(1+u)^2 + v^2}} \qquad
\]
has power series expansion
\begin{eqnarray*}
f(u,v) &=& (1+v^2(1+u)^{-2})^{-\frac{1}{2}}
~=~ \big(1+v^2-2uv^2+\mbox{O}(\rho^4)\big)^{-\frac{1}{2}} \\
&=& 1 - \frac{v^2}{2} + uv^2 +  \mbox{O}(\rho^4) \qquad \mbox{as}~ \rho \downarrow 0.
\end{eqnarray*}
The first term integrates to the area of the disk, $\pi \rho^2$.  The third term integrates to zero, being an odd function of $u$.  The second term, using symmetry in $u$ and $v$, integrates to
\[
\frac{1}{4} \int\!\!\!\int_{u^2+v^2 \le \rho^2} (u^2+v^2) \,du\,dv ~=~ 
\frac{1}{4}  \int_0^{2\pi} \!\! \int_0^\rho r^2\,r\,dr\,d\theta ~=~ \frac{1}{8} \pi\rho^4,
\]
so
\[
\beta ~=~ \pi \rho^2 \big(1 - \frac{\rho^2}{8} + \mbox{O}(\rho^4) \big),
\]
and the result follows.
\finpf

\begin{lem} \label{2}
For $\tau > 0$, if a random vector $X^\tau$ is uniformly distributed over the disk in $\R^2$, centered at the point $\tau e_1$ and of radius $(2\tau)^{3/2}$ , then, as $\tau \downarrow 0$,
\[
{\mathbb E}\frac{X^\tau_1}{|X^\tau|} ~=~ 1 - \tau + \mbox{\rm O}(\tau^2),
\]
and consequently
\[
{\mathbb E}\Big(X^\tau_1 + \frac{X^\tau_1}{|X^\tau|} - 1\Big) ~=~  \mbox{\rm O}(\tau^2).
\]
\end{lem}

\pf
The random vector $Y^\tau = \frac{1}{\tau}X^\tau$ is uniformly distributed over the disk in $\R^2$ of radius 
$\sqrt{8\tau}$ centered at the first unit vector $e_1$.  The result now follows from the previous lemma.
\finpf

With the assumptions of Lemma \ref{2}, symmetry implies
\[
{\mathbb E}\Big(X^\tau_2 + \frac{X^\tau_2}{|X^\tau|}\Big) ~=~  0.
\]
This observation, along with Lemma \ref{2}, proves ${\mathbb E}\big(\nabla f(X^\tau)\big) = \mbox{O}(\tau^2)$ in Example~\ref{counterexample}, as we claimed.

\section{Appendix:  The Goldstein modulus} \label{goldstein}
Following \cite{goldstein-modulus}, the {\em Goldstein modulus} of a locally Lipschitz function $f \colon \X \to \R$ at a point \mbox{$x \in \X$}, denoted $\Gamma f(x)$, is the infimum of those values $\e > 0$ such that some convex combination $y$ of subgradients in $\partial f\big(B_\e(x)\big)$ satisfies 
$|y| < \e$.  We say that the Goldstein modulus {\em grows linearly} at a point $\bar x \in \X$ if there exists a constant $\alpha > 0$ such that $\Gamma f(x) \ge \alpha|x-\bar x|$ for all $x$ near $\bar x$.

\begin{prop}
\mbox{}\!\!The Goldstein modulus for a locally Lipschitz function \mbox{$f \colon \X \to \R$} grows linearly at a point $\bar x \in \X$ if and only if $f$ has no Goldstein flat direction at $\bar x$.  In that case, the following two conditions also hold.
\begin{enumerate}
\item[{\rm (i)}]
The minimum norm of a Clarke subgradient grows linearly:  there exists a constant $\alpha > 0$ such that 
$|y| \ge \alpha|x-\bar x|$ for all $x$ near $\bar x$ and subgradients $y \in \partial f(x)$.
\item[{\rm (ii)}]
If $\bar x$ is a local minimizer for $f$, then $f$ grows quadratically at $\bar x$:  there exists a constant $\kappa > 0$ such that $f(x) - f(\bar x) \ge \kappa|x-\bar x|^2$ for all $x$ near $\bar x$.
\end{enumerate}
\end{prop}

\pf
If the function $f$ has a Goldstein flat direction $d$ at the point $\bar x$, then by definition there exist sequences of scalars $\tau_r > 0$ for $r=1,2,\ldots$ satisfying $\tau_r \downarrow 0$ as $r \to \infty$, and $\e_r = \mbox{o}(\tau_r)$, such that for each $r$, some convex combination $y_r$ of gradients of $f$ at points in the ball $B_{\e_r}(\bar x + \tau_r d)$ satisfies $|y_r| < \e_r$.  Since the Goldstein modulus therefore satisfies $\Gamma f(\bar x + \tau_r d) \le \e_r = \mbox{o}(\tau_r)$, it does not grow linearly at $\bar x$.

Conversely, suppose that the Goldstein modulus of the function $f$ does not grow linearly at the point 
$\bar x$, so there exists a sequence of scalars $\alpha_r > 0$ for $r=1,2,\ldots$ satisfying $\alpha_r \downarrow 0$ as $r \to \infty$, and a sequence of points $x_r \to \bar x$ in $\X$ such that, for all $r$ we have $\Gamma f(x_r) < \alpha_r |x_r - \bar x|$.  Define scalars $\tau_r = |x_r - \bar x|$.  We deduce the existence of a positive radius 
$\delta_r < \alpha_r \tau_r$ and a convex combination $y_r$ of subgradients of $f$ at points in the ball
$B_{\delta_r}(x_r)$ satisfying $|y_r| < \delta_r$. By increasing each $\delta_r$ slightly, we can suppose in addition that $y_r$ is a convex combination of gradients of $f$ at points in $B_{\delta_r}(x_r)$.  Clearly each $x_r$ is distinct from $\bar x$, so after taking a subsequence, we can suppose without loss of generality that the unit vectors 
$\tau_r^{-1}(x_r - \bar x)$ converge to some unit vector $d \in \X$, and hence
\[
|\bar x + \tau_r d - x_r| ~=~ \tau_r|d - \tau_r^{-1}(x_r - \bar x)| ~=~ \mbox{o}(\tau_r) \qquad \mbox{as}~ r \to \infty.
\]
The scalars $\e_r = \delta_r + |\bar x + \tau_r d - x_r|$ therefore also satisfy $\e_r = \mbox{o}(\tau_r)$ as $r \to \infty$.  The triangle inequality implies 
$B_{\delta_r}(x_r) \subset B_{\e_r} (\bar x + \tau_r d)$,
so each $y_r$ is a convex combination of gradients of $f$ at points in the ball
$B_{\e_r}(\bar x + \tau_r d)$ satisfying $|y_r| < \e_r$.  We deduce that $d$ is a Goldstein flat direction at $\bar x$, as required.

Suppose that the  Goldstein modulus of the function $f$ satisfies $\Gamma f(x) \ge \alpha|x-\bar x|$ for all points $x \in \X$ near $\bar x$.  The claim (i) follows immediately, because the Goldstein modulus satisfies
\[
\min_{y \in \partial f(x)} |y| ~\ge~ \Gamma f(x) ~\ge~ \alpha|x-\bar x|.
\]

Turning to the claim (ii), by assumption we have $0 \in \partial f(\bar x)$.  Consider any point 
$x \ne \bar x$ near $\bar x$.  Clearly the Goldstein modulus satisfies $\Gamma f(x) \le |x-\bar x|$. Fix any radius $\e$ in the interval $\big(0,\Gamma f(x)\big)$. The shortest convex combination $y$ of subgradients in $\partial f\big(B_\e(x)\big)$ satisfies $|y|>\e$.  Since $x$ is near $\bar x$, so is the point $x - \frac{\e}{|y|} y$,
so from \cite{gold_eps_stat} we see
\[
f(\bar x) ~\le~ f\Big(x - \frac{\e}{|y|} y\Big) ~\le~ f(x) - \e|y| ~<~ f(x) - \e^2.
\]
Taking the supremum over $\e$ shows
\[
f(x) \ge f(\bar x) + \big( \Gamma f(x) \big)^2 \ge \alpha^2|x-\bar x|^2,
\]
showing quadratic growth.
\finpf

\section{Appendix:  Proof of Example \ref{strongly-convex2}} \label{E}
We construct a semi-algebraic, strongly convex function $f \colon \R^2 \to \R$ such that 
\mbox{$0 \in \mbox{ri}\,\partial f(0)$} and yet $f$ is Goldstein flat at $0$.

Consider the convex function $g \colon \R^2 \to \R$ defined by
\[
g(x) ~=~ \max_{0 \le \theta \le 1,~\sigma=\pm 1} \left\{(1+\theta)\sigma x_2-\theta^2x_1-\theta^6\right\}.
\]
At the point $x = (0,0)$, the maximum is attained if and only if $\theta=0$, so
\[
\partial g(0,0) = \{0\} \times [-1,1].
\]
On the other hand, if $x_1 > 0$ and $0 < |x_2| < 6+2x_1$, then the maximum is attained uniquely by 
$\sigma = \mbox{sgn}(x_2)$ and with $\theta=\theta(x)$ being the unique solution of the equation  
$6\theta^5 + 2x_1 \theta = |x_2|$,
so in that case we have $\big(\nabla g(x)\big)_1 = -\theta^2(x)$.
The function $f = g + \frac{1}{2}|\cdot|^2$ is therefore strongly convex, with $(0,0)$ contained in the relative interior of the subdifferential
\[
\partial f (0,0) = \{0\} \times [-1,1],
\]
and, provided that $x_1 > 0$ and $0 < |x_2| < 6+2x_1$, we have $\big(\nabla f(x)\big)_1 = x_1-\theta^2(x)$.
Notice also that both $g$ and $f$ are even functions of $x_2$.

For small $\tau > 0$, suppose $x_1 = \tau + \mbox{O}(\tau^{3/2})$, and $x_2 = \mbox{O}(\tau^{3/2})$, and hence $|x_2|/x_1 = \mbox{O}(\tau^{1/2})$.  Since $6\theta^5 + 2x_1 \theta = |x_2|$, we deduce
\[
3(x_1^{-1/4}\theta)^5 + (x_1^{-1/4} \theta) ~=~ \frac{1}{2}x_1^{-5/4}|x_2| ~=~ \mbox{O}(\tau^{1/4})
\]
so via a power series expansion we see
\[
x_1^{-1/4} \theta ~=~ \frac{1}{2} x_1^{-5/4}|x_2| + \mbox{O}\big((x_1^{-5/4}|x_2|)^5\big).
\]
We deduce
\[
\theta(x) = \frac{|x_2|}{2x_1} + \mbox{O}(\tau^{3/2})
\]
so
\[
\theta^2(x) ~=~ \frac{x^2_2}{4x^2_1} + \mbox{O}(\tau^2) ~=~ \frac{x^2_2}{4\big(\tau(1+\mbox{O}(\tau^{1/2})\big)^2} + \mbox{O}(\tau^2) ~=~ \frac{x^2_2}{4\tau^2} + \mbox{o}({\tau}).
\]
Now consider a random vector $X^\tau$ uniformly distributed on a ball of radius $4\tau^{3/2}$ centered at the point $(\tau,0)$.  A standard calculation in polar coordinates shows
\[
{\mathbb E}\big((X^\tau_2)^2\big) ~=~ \frac{1}{4} (4\tau^{3/2})^2 = 4\tau^3,
\]
so ${\mathbb E}\big(\theta^2(X^\tau)\big) = \tau + o(\tau)$,
and hence 
${\mathbb E}\big(\nabla f(X^\tau)\big)_1 = o(\tau)$.
By symmetry, we have ${\mathbb E}\big(\nabla f(X^\tau)\big)_2 = 0$,
so we conclude ${\mathbb E}\big(\nabla f(X^\tau)\big) = o(\tau)$.
This shows that $(1,0)$ is a Goldstein flat direction for $f$ at the minimizer $(0,0)$.

\section{Appendix:  Proof of Proposition \ref{plq}} \label{D}
Without loss of generality, consider a continuous, piecewise linear-quadratic function $f \colon \Rn \to \R$ and suppose that the point $0 \in \Rn$ is a strict local minimizer.  We can suppose that there exist polyhedral cones $K_i \subset \Rn$ for $i$ in some finite index set $I$, along with corresponding vectors $b_i \in \Rn$ and symmetric $n$-by-$n$ matrices $A_i$, such that $\cup_i K_i = \Rn$ and the quadratic functions $f_i \colon \Rn \to \R$ defined by
\[
f_i(x) ~=~ \frac{1}{2} x^T A_i x + b_i^T x \qquad (x \in \Rn)
\]
satisfy $f(x)= f_i(x)$ for all points $x \in K_i$.

Suppose that there exists a Goldstein flat direction $d \in \Rn$ for $f$ at $0$.  By definition, there exist sequences of scalars 
$\tau_r > 0$ for $r=1,2,\ldots$ satisfying $\tau_r \downarrow 0$ as $r \to \infty$, and $\e_r = \mbox{o}(\tau_r)$, such that for each $r$, some convex combination $y_r$ of gradients of $f$ at points in the ball $B_{\e_r}(\tau_r d)$ satisfies $|y_r| < \e_r$. Define a nonempty set
\[
\hat I = \{i \in I : d \in K_i \}. 
\] 
For all $i \not\in \hat I$, we have
\[
K_i \cap B_{\e_r}(\tau_r d) ~=~ K_i \cap \tau_r B_{\e_r/\tau_r}(d) ~=~ 
\tau_r\big(K_i \cap B_{\e_r/\tau_r}(d)\big) ~=~ 
\emptyset \qquad \mbox{for all large}~r.
\]
On the other hand, for all $i \in \hat I$ we have
\[
f(\tau d) ~=~ \frac{1}{2} (d^T A_i d)\tau^2  +  (b_i^T d)\tau \qquad \mbox{for all}~ \tau \ge 0.
\] 
Consequently, there must exist scalars $\gamma$ and $\mu$ such that $b_i^T d = \gamma$ and 
$d^T A_i d = \mu$ for all $i \in \hat I$, and furthermore
\[
f(\tau d) ~=~ \frac{1}{2} \mu\tau^2  +  \gamma\tau \qquad \mbox{for all}~ \tau \ge 0.
\]
Since $0$ is a strict local minimizer, either $\gamma > 0$, or $\gamma = 0$ and $\mu > 0$.

By assumption, for all $r=1,2,\ldots$, we have
\begin{eqnarray*}
y_r 
&\in& 
\mbox{conv} \bigcup_{i \in \hat I} \nabla f_i\big(B_{\e_r}(\tau_r d)\big)
~=~ 
\mbox{conv} \bigcup_{i \in \hat I} (b_i + A_i\big(B_{\e_r}(\tau_r d)\big)) \\
&=& 
\mbox{conv} \bigcup_{i \in \hat I} (b_i + \tau_r A_i d + \epsilon_r A_i B),
\end{eqnarray*}
where $B$ denotes the closed unit ball in $\Rn$.  Letting $r \to \infty$ shows
\[
0 ~\in~ \mbox{conv}\{b_i : i \in \hat I\}
\]
which implies $\gamma = 0$.  Consequently we must have $\mu > 0$.  We deduce
\[
d^T y_r ~\in~ \mbox{conv} \bigcup_{i \in \hat I} (\tau_r \mu + \epsilon_r d^TA_i B),
\]
so $|d^T y_r - \tau_r \mu| = O(\e_r) = o(\tau_r)$, leading to the contradiction 
\[
0 = \lim_{r \to \infty} \frac{d^T y_r}{\tau_r} = \mu > 0.
\]
This completes the proof.

\section{Appendix:  Proof of Example \ref{fully-amenable1}} \label{F}
Consider the convex function $f \colon \R^2 \to \R$ defined by
\[
f(x,y) ~=~ \max\{x^2+y^2,2|y|\Big\} \qquad \mbox{for}~ (x,y) \in \R^2.
\]
The function grows quadratically at its global minimizer $(0,0)$, which is in fact tilt stable \cite{tilt_orig}:  specifically, for all small vectors $(u,v) \in \R^2$, the tilted function
\[
(x,y) ~\mapsto~ f(x,y) - 2ux - 2vy
\]
 has a unique minimizer, whose dependence on $(u,v)$ is locally Lipschitz.

To verify this claim, we first note that the tilted function has compact lower level sets, so has a minimizer.  Define a function \mbox{$\alpha \colon (-1,1) \to \R$} by
\[
\alpha(x) = 1-\sqrt{1-x^2} \qquad (-1 < x < 1).
\]
Notice that $\alpha$ is locally Lipschitz.
Provided that $|y| < 1$ we have
\[
f(x,y) - 2vy ~=~ 
\left\{
\begin{array}{ll}
x^2+y^2 - 2vy	& \mbox{if}~ |y| \le \alpha(x) \\
2|y| - 2vy		& \mbox{otherwise}.
\end{array}
\right.
\]
For any fixed $x \in (-1,1)$, this function is a piecewise-quadratic convex function of the variable 
$y \in \R$.  For $v \ge 0$, the unique minimizer is $y = \min\{v,\alpha(x)\}$, so by symmetry, in general the minimizer is
\[
y ~=~ \mbox{sgn}(v) \min\{ |v|,\alpha(x) \}.
\]
Notice that $y$ is a locally Lipschitz function of $x$ and $v$.
Substituting into the tilted function gives the value
\[
\left\{
\begin{array}{ll}
x^2-v^2 -2ux 
& \mbox{if}~ x^2 \ge 2|v|-v^2 \\
2(1-|v|)\alpha(x) - 2ux 
& \mbox{otherwise}.
\end{array}
\right.
\]
We assume $|v| < 1$.  In that case, this function of $x$ is strongly convex, so it has a unique minimizer, given by the continuous function
\[
x ~=~
\left\{
\begin{array}{ll}
u & \mbox{if}~ u^2 \ge 2|v| - v^2  \\
p(u,v) & \mbox{otherwise}.
\end{array}
\right.
\]
where $p(u,v)$ is the unique solution $x$ of the equation
\[
\alpha'(x) = \frac{u}{1-|v|},
\]
namely
\[
x = p(u,v) = \frac{u}{\sqrt{u^2 + (1-|v|)^2}}.
\]
The dependence of $x$ on $(u,v)$ is locally Lipschitz, and hence so is that of~$y$.

To summarize, the function $f$ is convex, semi-algebraic, piecewise ${\mathcal C}^{(2)}$, and fully amenable, and grows quadratically around the tilt-stable global minimizer at $(0,0)$.  Nonetheless, $f$ is Goldstein flat at $(0,0)$.  Indeed, for $\tau > 0$, if a random vector $X^\tau$ is uniformly distributed over the disk in $\R^2$, centered at the point $\tau e_1$ and of radius $\tau^{3/2}$, then, as $\tau \downarrow 0$,
\[
{\mathbb E}\big(\nabla f(X^\tau)\big) ~=~ \mbox{O}(\tau^{3/2}).
\]
To see this, note by symmetry, we have
\[
{\mathbb E}\big(\nabla f(X^\tau)_2\big) ~=~ 0.
\]
The area of the disk is $\pi \tau^3$.  Within the disk, define a set
\[
S ~=~ \{ (x,y) : x^2+y^2 > 2|y| \}.
\]
This set has area $\mbox{O}(\tau^{7/2})$, so the integral over $S$ of $\nabla f(x)_1 = 2x_1$ is  $\mbox{O}(\tau^{9/2})$.  On the other hand, the corresponding integral over the complement of $S$ is zero, so the result follows.

\end{document}